\documentclass{IEEEtran4PSCC}

\usepackage{etoolbox}
\makeatletter
\patchcmd{\@makecaption}
{\scshape}
{}
{}
{}
\makeatletter
\patchcmd{\@makecaption}
{\\}
{.\ }
{}
{}
\makeatother
\def\tablename{Table}

\ifCLASSINFOpdf
   \usepackage[pdftex]{graphicx}
\else
   \usepackage[dvips]{graphicx}
\fi
\usepackage[cmex10]{amsmath}
\usepackage{xcolor}
\makeatletter
\let\old@ps@headings\ps@headings
\let\old@ps@IEEEtitlepagestyle\ps@IEEEtitlepagestyle
\def\psccfooter#1{%
    \def\ps@headings{%
        \old@ps@headings%
        \def\@oddfoot{\strut\hfill#1\hfill\strut}%
        \def\@evenfoot{\strut\hfill#1\hfill\strut}%
    }%
    \def\ps@IEEEtitlepagestyle{%
        \old@ps@IEEEtitlepagestyle%
        \def\@oddfoot{\strut\hfill#1\hfill\strut}%
        \def\@evenfoot{\strut\hfill#1\hfill\strut}%
    }%
    \ps@headings%
}
\makeatother

\psccfooter{%
        \parbox{\textwidth}{\hrulefill \\ \small{22nd Power Systems Computation Conference} \hfill \begin{minipage}{0.2\textwidth}\centering \vspace*{4pt} \includegraphics[scale=0.06]{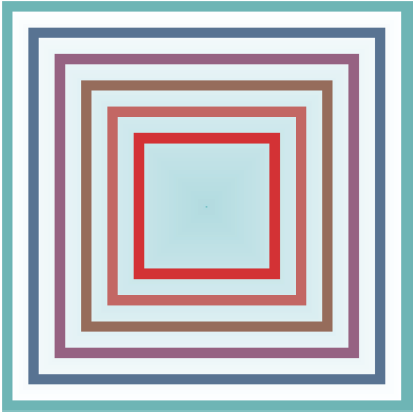}\\\small{PSCC 2022} \end{minipage} \hfill \small{Porto, Portugal --- June 27 -- July 1, 2022}}%
}

\usepackage{cite,graphics}

\renewcommand{\baselinestretch}{1}
\usepackage{amsmath,amsfonts,amssymb,mathtools,amsthm}
\providecommand{\n}[1]{\lVert#1\rVert}

\providecommand{\tx}[1]{\text{\upshape{#1}}}

\usepackage[printonlyused]{acronym}
\acrodef{hac}[HAC]{hybrid angle control}
\acrodef{coi}[COI]{center-of-inertia}
\acrodef{ib}[IB]{infinite bus}
\acrodef{sm}[SM]{synchronous machines}
\acrodef{wrt}[w.r.t.]{with respect to}
\acrodef{agas}[AGAS]{almost global asymptotic stability}
\acrodef{lhs}[LHS]{left-hand side}  
\acrodef{rhs}[RHS]{right-hand side}  
\acrodef{rocof}[RoCoF]{rate of change of frequency}  
\acrodef{avr}[AVR]{automatic voltage regulator}  
\acrodef{pss}[PSS]{power system stabilizer}
\acrodef{gfc}[GFC]{grid-forming converter}
\acrodef{vsm}[VSM]{virtual synchronous machine}
\acrodef{dvoc}[dVOC]{dispatchable virtual oscillator control}

\newtheorem{remark}{Remark}

\usepackage{academicons,fontawesome}
\begin{document}
%
\title{System-Level Performance and Robustness of the Grid-Forming Hybrid Angle Control}

\author{
\IEEEauthorblockN{Ali Tayyebi$^{*\dag}$, Alan Magdaleno$^{\dag}$, Denis Vettoretti$^{*}$, Meng Chen$^{\ddag}$, \\Eduardo Prieto-Araujo$^{\star}$, Adolfo Anta$^{*}$, and Florian Dörfler$^{\dag}$}
\IEEEauthorblockA{$^{*}$Electric Energy Systems, Austrian Institute of Technology, Vienna, Austria.\\
$^{\dag}$Automatic Control Laboratory (IfA), Swiss Federal Institute of Technology (ETH), Zürich, Switzerland.\\
$^{\ddag}$AAU Energy, Aalborg University, Aalborg, Denmark\\
$^\star$CITCEA, Universitat Politècnica de Catalunya (UPC), Barcelona, Spain.\\
Corresponding author's e-mail address: ali.tayyebi-khameneh@ait.ac.at
}
}


\maketitle

\begin{abstract}
This paper investigates the implementation and application of the multi-variable grid-forming \ac{hac} for high-power converters in transmission grids. We explore the system-level performance and robustness of \ac{hac} concept in contrast to other grid-forming schemes i.e., power-frequency droop and matching controls. Our findings suggests that, similar to the ac-based droop control, \ac{hac} enhances the small-signal frequency stability in low-inertia power grids, and akin to the dc-based matching control, \ac{hac} exhibits robustness when accounting for the practical limits of the converter systems. Thus, \ac{hac} combines the aforementioned complementary advantageous. Furthermore, we show how retuning certain control parameters of the grid-forming controls improves the frequency performance. Last, as separate contributions, we introduce an alternative control augmentation that enhances the robustness and provide theoretical guidelines on extending the stability certificates of \ac{hac} to multi-converter systems.     
\end{abstract}

%



\section{Introduction}
The ambitious targets that are set to globally reduce the carbon footprint require revolutionizing the foundations of legacy power systems. In other words, the \ac{sm}-based energy generation from fossil fuels must be replaced with green and sustainable energy sources. The majority of clean energy sources interface the power grids via dc-ac power converters. Although, in contrast to the bulk \ac{sm}s, power converters are fast, modular, and highly controllable, they are subject to volatile energy resources  and lack the necessary robustness and maturity to ensure adequate and reliable power delivery. It is envisioned that the advanced control architectures for power converters can possibly address the aforementioned concerns \cite{MDHHV18,kroposki2017achieving,RWLLE21,8070502}. 

The so-called \emph{grid-following} converter controls that exploit an explicit synchronizing mechanism are widely utilized \cite{li2021rethinking,8586162,TDKZH18}. However, the grid-following converters exhibit robustness and stability issues in the converter-dominated grids that are highlighted by a significant reduction of rotational inertia i.e., \emph{low-inertia} grids; see \cite{MOVAH19,CTGAKF19,8973657} among others. Subsequently, the concept of \emph{\ac{gfc}} is introduced that provides fast and robust frequency and voltage regulation to address the stability challenges associated with the low-inertia systems \cite{TGAKD20,UNSW20}. 

Several grid-forming control techniques have been proposed in recent years. Restricting the focus to converter frequency definition under these controls, one can highlight fundamental differences. The vast majority of grid-forming techniques e.g., droop control, \ac{vsm}, and oscillator-based schemes shape the converter frequency based on the ac quantities such as current, voltage, and power flows e.g., \cite{JMAFD:15,SDJD17,CZB21,ZW11,LKSW21,CDA93,CGBF19,DSF15,SDB13,SJB21,rowe2012arctan}. Recently, a class of controllers have been proposed that define the converter frequency in proportion to a linear/nonlinear dc voltage feedback \cite{AJD18,CGD17,cvetkovic_modeling_2015,huang2017virtual,WL20,AF20,JD20}. Most recently, the emerging multi-variable control trend is to combine dc and ac information in designing the \ac{gfc} frequency \cite{TAD20,TAC2020,GRL20,SG21,chen2021generalized}.

In this work, we focus on grid-forming \acf{hac} that blends linear dc and nonlinear ac feedback for defining converter frequency \cite{TAD20,TAC2020}. We provide details on the implementation of \ac{hac}, highlight its performance in contrast to other grid-forming controls via various simulation case studies, and propose a complementary control augmentation for enhancing the robustness. 
In addition, a high-level stability analysis is presented that provides guidelines on extending the stability certificates of \ac{hac} (as in \cite{TAC2020}) to an interconnected system of converters. Finally, as a separate contribution in its own right, we show that an appropriate tuning of converter control parameters can globally enhance frequency stability across the low-inertia power system.     

The remainder of this paper is structured as it follows: Section \ref{sec:modelling} presents our modeling approach for low-inertia power grids, Section \ref{sec:control} describes three different grid-forming controls and highlights their dominant features. Furthermore, it provides practical details on the implementation of \ac{hac}. Section \ref{sec:simulation} presents several comparative case studies. Last, the Appendix introduces a variant of \ac{hac} with enhanced robustness and includes a high-level stability investigation for multi-converter systems under \ac{hac}.

\begin{figure*}
	\centering	
	{\includegraphics[trim=2.6mm 1.65mm 1.1mm 1.1mm ,clip,width=1.5\columnwidth]{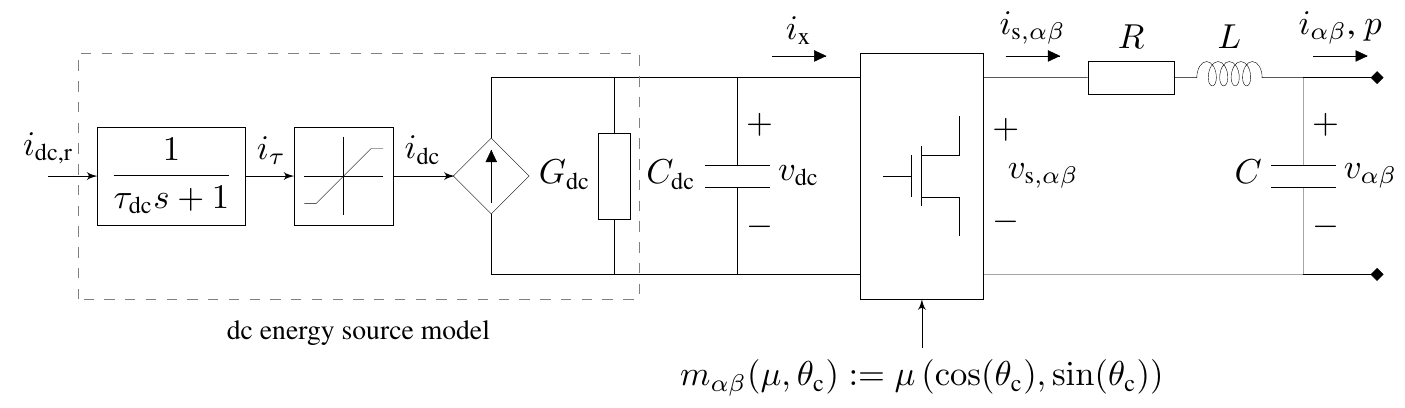}}
	\caption{DC-AC power converter model schematic with ac quantities represented in stationary $\alpha\beta$-coordinates \cite{TGAKD20,YI10}.}\label{fig:converter}
\end{figure*}
\section{Modeling Description}\label{sec:modelling}
{In this paper, we consider a transmission power grid that is fed by a mixture of \ac{sm} and converter-based generation units}. For the sake of completeness, we briefly review the \ac{sm}, network, and converter modeling; the reader is referred to \cite{TGAKD20} for further elaborations.

\subsection{Synchronous Machine}
In this work, we consider a detailed nonlinear \ac{sm} model that incorporates second-order mechanical dynamics (associated with angle and frequency) and sixth-order electrical dynamics (associated with the stator, damper, and field winding fluxes). Furthermore, the \ac{sm} is combined with a first-order dynamic turbine model. The reader is referred to \cite[Section II-B]{TGAKD20}\cite{kundur1994power} for a detailed description of the SM modeling. It is noteworthy that we include the standard \ac{sm} control mechanisms such as \ac{avr}, \ac{pss}, and turbine governor control; see \cite{TGAKD20} for details. Furthermore, the governor gain in \cite[Equation 5]{TGAKD20} is exploited to realize coordinated load-sharing with other generation units across the system; see \cite[Appendix]{TGAKD20} for a tuning criteria recommendation. 
\subsection{Power Network}
The generation units {are considered to} interface the transmission lines via identical medium-to-high voltage (MV/HV) transformers. We consider dynamic transmission lines represented with standard $\pi$-section models \cite{kundur1994power,TGAKD20}. Furthermore, the network loads are modeled by constant impedances. {We emphasize that the load models are voltage-dependent and restrict our focus to the active power loads that primarily influence the frequency dynamics in transmission grids}. Note that due to the fast timescales of the \ac{gfc}s and potential adverse interactions with line dynamics, the quasi-steady-state algebraic network model assumption is not valid \cite[Subsection III-A]{CGBF19}, hence it is necessary to consider dynamic line model{s for low-inertia power systems} \cite{MOVAH19,GCBD19}. 

\subsection{High-Power DC-AC Converter}
The power converter model in this work includes an explicit first-order representation of the dc energy source dynamics, dc-link capacitance, power-preserving (i.e., lossless) average model of the {two-level dc-ac conversion stage}, and the LC output filter. 
Furthermore, we consider the current limits of the dc energy source that in practice corresponds to the current limits of the PV, battery, or wind generator systems. {It is worth mentioning that in this work the dc source model represents an aggregation of several energy sources that supply a modular high-power converter; see Section \ref{sec:simulation} for further details}. Figure \ref{fig:converter} present the schematic of the converter model; see \cite[Subsection II-A]{TGAKD20} for further details and the differential equations corresponding to the Figure \ref{fig:converter}.

The converter model illustrated in Figure \ref{fig:converter} provides three control degrees of freedom, namely: 1) the dc source current reference $i_\tx{dc,r}$, 2) the modulation signal magnitude i.e., $\mu=\n{m_{\alpha\beta}(\mu,\theta_\tx{c})}$, and the modulation signal angle i.e., $\theta_\text{c}=\angle m_{\alpha\beta}(\mu,\theta_\tx{c})$. In Section \ref{sec:control}, we elaborate on the definition of these control inputs under three different grid-forming control strategies.
\section{Grid-Forming Frequency Controls}\label{sec:control}
The vast literature on grid-forming control schemes can be distinctly assigned into three categories: 
\begin{enumerate}
\item {ac-based} control techniques that define the converter frequency based on ac measurements e.g., see  \cite{JMAFD:15,SDJD17,CZB21,ZW11,LKSW21,CGBF19,DSF15,CDA93,SDB13,SJB21,rowe2012arctan},
\item {dc-based} model-matching inspired techniques that relate the converter frequency to the dc measurements e.g., see \cite{AJD18,CGD17,cvetkovic_modeling_2015,huang2017virtual,WL20,AF20,JD20}, and
\item {multi-variable hybrid} control structures that exploit both dc and ac measurements for a grid-forming frequency synthesis e.g., see \cite{TAD20,TAC2020,GRL20,SG21,chen2021generalized}.
\end{enumerate}
In what follows, we review a candidate from each control category and {highlight their dominant features}.


\subsection{AC-Based Power-Frequency Droop Control}
The baseline power-frequency droop control that is inspired by the frequency droop behavior of the \ac{sm}, serves as a powerful yet simple control solution for grid-forming converter applications \cite{CDA93,SDB13,SJB21,rowe2012arctan}. The power-frequency droop control in its simplest form is described by
\begin{equation}\label{eq:droop}
\dot{\theta}_\tx{c}=\omega_\tx{c}:=\omega_0+d_{p-\omega}(p_\tx{r}-p),
\end{equation}
where $\theta_\tx{c}$ denotes the phase angle of the converter modulation signal $m_{\alpha\beta}(\mu,\theta_\tx{c})$ in Figure \ref{fig:converter}, $\omega_\tx{c}$ denotes the converter angular frequency, $\omega_0$ denotes the nominal system frequency, $d_{p-\omega}$ is the power-frequency droop gain, $p_\tx{r}$ is the control reference for the active power flowing out of the converter's ac-side terminal that is defined by $p:=v^\top_{\alpha\beta}i^{}_{\alpha\beta}$. The ac-based nature of the droop control can be observed in \eqref{eq:droop} that highlights the frequency dependency on the ac active power feedback. It is worth mentioning that low-pass filtering the active power feedback in \eqref{eq:droop} is a common practice to safeguard the converter control against measurement imperfections and high-frequency harmonics \cite[Introduction]{JRSD17}. Section \ref{sec:simulation} presents a case study that unveils the strong influence of such low-pass filtering in reshaping the system-level post-contingency frequency evolution. 
\begin{remark}\textup{(Enhanced small-signal frequency stability)}\label{rem:ac}\\
{Recent explorations \cite{TGAKD20,CTGAKF19,TDKZH18} uncover that} the ac-based grid-forming control strategies e.g., droop control, \ac{vsm}, and \ac{dvoc} clearly improve the small-signal frequency stability (in terms of the \ac{rocof} and nadir performance metrics) of the low-inertia system compared to the all-SMs conventional system. This is primarily underpinned by the fast response timescale of the \ac{gfc}s that enables fast frequency regulation. Furthermore, it has been observed that the ac-based grid-forming schemes exhibit slightly better frequency response compared to the dc-based counterparts e.g., the matching control (presented in the next Subsection); see \cite[Subsection IV-C]{TGAKD20} for details. This behavior is due to the fact that following a network contingency/disturbance, the ac feedback (e.g., $p$ in \eqref{eq:droop}) quickly reflects the grid conditions in the converter frequency definition and enables fast frequency regulation.
\end{remark}  
\subsection{DC-Based Model-Matching Control}
Inspired by the structural similarities of the \ac{sm} and converter dynamical models, a family of model-matching control techniques have been proposed \cite{AJD18,CGD17,cvetkovic_modeling_2015,huang2017virtual,WL20,AF20,JD20}. These control structures establish strong duality between the converter and \ac{sm} dynamical structure; see \cite{AF20} for a detailed derivation. Matching control in its original form is described by \cite{AJD18}
\begin{equation}\label{eq:matching}
\dot{\theta}_\text{c}=\omega_\tx{c}:=\eta v_\tx{dc},
\end{equation} 
where $\eta:=\omega_0/v_\tx{dc,r}$ in which $v_\tx{dc,r}$ denotes the reference for the dc voltage $v_\tx{dc}$. Note that \eqref{eq:matching} highlights the dc-based nature of the matching strategy. Similar to droop control, it might be necessary to low-pass filter the dc feedback in \eqref{eq:matching} to safeguard the angle dynamics against potential dc-link voltage ripples.
\begin{remark}\textup{(Enhanced robustness)}\label{rem:dc}\\
Recall that the \ac{sm}s can be operated with either flexible or constant mechanical input torque, and regardless of this degree of freedom, the \ac{sm}s achieve  robust synchronization with the power grid (thanks to their inherent self-synchronizing feature \cite{BSEO17,cvetkovic_modeling_2015,AJD18,AF20}).  	
Under the strong converter-\ac{sm} duality induced by matching control, the converter can be operated with either constant or flexible input dc current that is injected by the dc energy source in Figure \ref{fig:converter}. This property of the matching control comes strongly into the picture when considering the safety current constraints of the dc energy source; see the saturation function in Figure \ref{fig:converter}. Recent works \cite{TGAKD20,SC21,GRL20} formally and numerically demonstrate that the  matching controls (unlike the ac-based schemes) preserve the closed-loop stability under active dc current constraint. This behavior can be perceived as a built-in mode-switching feature that simultaneously changes the converter operation mode from voltage (i.e., grid-forming) to current (i.e., grid-following) source when the dc constraint is activated; e.g., see \cite[Figure 15]{TGAKD20}. To our knowledge this is the most dominant feature of the dc-based controls i.e., superior robustness.
\end{remark}  
\subsection{Hybrid Angle Control}\label{subsec:hac}
The hybrid ac/dc grid-forming control architecture appears as a natural extension and the promising solution for combining the aforementioned complementary benefits of the ac/dc-based control schemes \cite{TAD20,TAC2020,GRL20,SG21}; see Remarks \ref{rem:ac} and \ref{rem:dc}. In what follows, we restrict our focus to the recently proposed \ac{hac} that blends linear dc voltage and nonlinear ac angle error terms for defining the converter frequency \cite{TAD20,TAC2020}. The grid-forming \ac{hac} is described by
\begin{equation}\label{eq:hac}
\dot{\theta}_\tx{c}=\omega_\tx{c}:=\omega_0+\underbracket{\gamma_\tx{dc}\left(v_\tx{dc}-v_\tx{dc,r}\right)}_{\tx{dc model-matching}}-\underbracket{\gamma_\tx{ac}\sin\left(\dfrac{\delta-\delta_\tx{r}}{2}\right)}_\tx{ac angle synchronization},
\end{equation} 

where $\gamma_\tx{dc}$ and $\gamma_\tx{ac}$ respectively denote the dc and ac control gains, and $\delta:=\angle v_{\text{s},\alpha\beta}-\angle v_{\alpha\beta}=\theta_\tx{c}-\theta_v$ denotes the phase angle difference (i.e., the relative angle) of the voltage before the filter and the output voltage in Figure \ref{fig:converter}, and $\delta_\tx{r}$ denotes the control reference for $\delta$. 
Note that the dc term in \eqref{eq:hac} is identical to the matching control variant proposed in \cite{CGD17} and the nonlinear ac term in \eqref{eq:hac} resembles the Kuramoto-like angle synchronizing term associated with the classic droop control \eqref{eq:droop} under certain assumption{s} \cite[Remark 3]{TAC2020}\cite{SDB13}. The practical implementation of \eqref{eq:hac} will be explored in the next section. The reader is referred to \cite{TAD20,TAC2020} for details on the design and properties of the \ac{hac}. {However, the following presents a brief summary:} 
\begin{enumerate}
\item \ac{hac} \eqref{eq:hac} incorporates an inherent dc-ac power-balancing behavior i.e., assuming $\omega_\tx{c}\to\omega_0$ then $\gamma_\tx{dc}\left(v_\tx{dc}-v_\tx{dc,r}\right)-\gamma_\tx{ac}\sin\left({\delta-\delta_\tx{r}}/{2}\right)\to0$. For instance, this means that if $v_\tx{dc}>v_\tx{dc,r}$ then $\delta>\delta_\tx{r}$ that allows for increased power injection into the ac-side and subsequently stabilizes the dc voltage.
\item The dc gain $\gamma_\tx{dc}$ predominantly reinforces the frequency dependency on the dc dynamics. {It is expected that a nonzero dc gain enables \ac{hac} to exhibit the robustness of dc-based controls (see Remark \ref{rem:dc} for details and Section \ref{sec:simulation} for a numerical justification of this hypothesis.)}
\item The ac gain $\gamma_\tx{ac}$ strongly influences the timescale associated with the ac-side power flows. {Therefore, it resembles the influence of the droop gain on the performance of the droop control.} 
\end{enumerate}

\subsection{DC and AC Voltage Control {Schemes}}
The grid-forming droop, matching and hybrid angle controls \eqref{eq:droop}-\eqref{eq:hac} define the phase angle that enters the converter modulation signal $m_{\alpha\beta}(\mu,\theta_\text{c})$ in Figure \ref{fig:converter}. It remains to close the loop by assigning the remaining control inputs i.e., the dc energy source reference current $i_\tx{dc,r}$ and modulation signal magnitude $\mu$. For the sake of fairness in the forthcoming comparative investigation, we consider identical complementary dc and ac voltage controls (that respectively define $i_\tx{dc,r}$ and $\mu$) for the aforementioned grid-forming strategies. 

Concerning the dc voltage control, we adopt the scheme that is proposed in \cite{TGAKD20} (see \cite{AF20} for a similar approach) i.e.,
\begin{align}\label{eq:i_dcr}
{i}_{{\text{\upshape{dc,r}}}}:=\underbracket{\kappa_{{\text{\upshape{dc}}}}\left(v_{{\text{\upshape{dc,r}}}}-{v_{{\text{\upshape{dc}}}}}\right)}_{\text{proportional control}}+\underbracket{\frac{p_\tx{r}}{v_{{\text{\upshape{dc,r}}}}}+\left(G_{{\text{\upshape{dc}}}}{v_{{\text{\upshape{dc}}}}}+\frac{v_{{\text{\upshape{dc}}}}i_{\text{\upshape{x}}}-{{p}}}{v_{{\text{\upshape{dc,r}}}}}\right)}_{\text{power injection and loss feedforward}}
\end{align}
where $k_\tx{dc}$ denotes the proportional dc voltage control gain and $i_\tx{x}:=m^\top_{\alpha\beta}i^{}_{\alpha\beta}$ denotes the net dc current injection to the ac-side; see Figure \ref{fig:converter}. Note that the power injection and loss compensation terms\footnote{The compensation scheme in \eqref{eq:i_dcr} is conventionally known as feedforward control, however due to the presence of the state-dependent quantities (such as $v_\tx{dc}$, $i_\tx{x}$, and $p$) it technically represents an algebraic state feedback.}  in \eqref{eq:i_dcr} are not necessary but improve the dynamic response and power set-point tracking \cite{TGAKD20,AF20}.

Next, inspired by the \ac{avr} mechanism of the \ac{sm}s, we augment the grid-forming strategies with a proportional-integral (PI) ac voltage magnitude control \cite{TGAKD20,SKBH21} that is (with slight abuse of the notation) given by 
\begin{equation}\label{eq:ac voltage}
\mu:=\kappa_{\text{p}}\left(v_\tx{r}-\n{v_{\alpha\beta}(t)}\right)+\kappa_{\text{i}}\int_{0}^{t}\left(v_\tx{r}-\n{v_{\alpha\beta}(s)}\right) \text{d}s,
\end{equation}
where $\kappa_\tx{p}$ and $\kappa_\tx{i}$ respectively denote the proportional and integral ac voltage control gains, and $v_\tx{r}$ denotes the reference ac voltage amplitude. We remark that the robust ac voltage regulation under \eqref{eq:ac voltage} is favorable for a \ac{gfc} (particularly in islanded/microgrid applications). However, one can alternatively consider a classic reactive power-voltage droop control as in \cite{MOVAH19}. Subsequently, the converter modulation signal in Figure \ref{fig:converter} is defined based the angle in \eqref{eq:droop}-\eqref{eq:hac} and the magnitude prescribed by \eqref{eq:ac voltage}. Last, in case studies presented in the Section \ref{sec:simulation} we identically tune the dc and ac voltage controls \eqref{eq:i_dcr} and \eqref{eq:ac voltage} for all grid-forming strategies.  

\subsection{Hybrid Angle Control Implementation}
In this subsection, we elaborate on the implementation of the grid-forming \ac{hac}. Before embarking upon these discussions, we remark that the implementation of droop control \eqref{eq:droop}, matching control \eqref{eq:matching}, and voltage controls \eqref{eq:i_dcr} and \eqref{eq:ac voltage} is previously addressed; the reader is referred to \cite{TGAKD20,TGA19,AF20,CDA93,CGD17} for details.

Concerning the implementation of \ac{hac} \eqref{eq:hac}, previous works \cite{TAD20,TAC2020} establish the theoretical foundations. However, the implementation presented in \cite{TAC2020,TAD20} is not straightforward, thus in what follows, we make the reported ideas clear. 

\begin{figure*}
	\centering	
	{\includegraphics[trim=17mm 0mm 0mm 0mm ,clip,width=1.5\columnwidth]{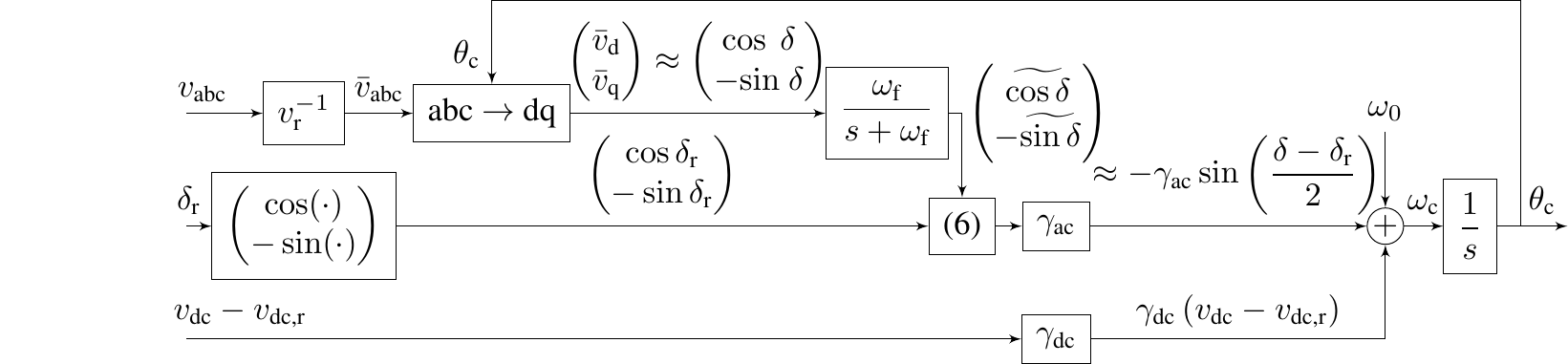}}
	\caption{Block diagram of the \ac{hac} \eqref{eq:hac} based on the dc voltage and the output ac voltage measurements in Figure \ref{fig:converter}.}\label{fig:hac}
\end{figure*}

To begin with, the dc component of the \ac{hac} \eqref{eq:hac} is easily constructed based on the dc voltage measurement in Figure \ref{fig:converter}; see the the dc feedback control in Figure \ref{fig:hac}.
Prior to describing the implementation of the angle synchronizing term in \eqref{eq:hac}, note that $\delta=\theta_\tx{c}-\theta_v$, thus one has to implicitly/explicitly derive the angle information from the output ac voltage measurement. Assume that the ac voltage control \eqref{eq:ac voltage} is sufficiently fast such that $\n{v_{\alpha\beta}}\approx v_\tx{r}$ (i.e., $\n{v_\tx{abc}}\approx v_\tx{r}$ that follows from the magnitude-preserving Clarke transformation \cite{YI10,TAC2020}). 

As it is illustrated in Figure \ref{fig:hac}, $\bar{v}_\tx{abc}:=v_\tx{abc}/v_\tx{r}$ denotes the normalized three-phase ac output voltage (note that $\bar{v}_\tx{abc}$ represents a unity phasor in polar coordinates that rotates with angle $\theta_v$). Next, we transform this quantity to a dq-coordinates that is aligned with the converter modulation angle $\theta_\tx{c}$. Thus, the image of $\bar{v}_\tx{abc}$ in the dq-rectangular coordinates (i.e., the unit vector rotating with the angle $\theta_v-\theta_\tx{c}{=-\delta}$) implicitly contains the relative angle information i.e., 
\begin{equation*}
{\left(\bar{v}_\tx{d},\bar{v}_\tx{q}\right)= \left(\cos(\theta_v-\theta_\tx{c}),\sin(\theta_v-\theta_\tx{c})\right) =\left(\cos\delta,-\sin\delta\right).}
\end{equation*}
Similar to low-pass filtering the active power feedback for the droop control \eqref{eq:droop}, we apply a first-order low-pass filter to $\left(\cos\delta,-\sin\delta\right)$ where $\omega_\tx{f}$ denotes the cutoff frequency. Subsequently, the filter output represents an approximation of the $\left(\cos\delta,-\sin\delta\right)$ that is $\left(\widetilde{\cos\delta},-\widetilde{\sin\delta}\right)$\footnote{Note that the LPF can by alternatively applied to $v_\tx{abc}$ in Figure \ref{fig:hac}.}. 
On the other hand, we process a given relative angle reference $\delta_\tx{r}$ that is consistent with desired power flows by the trigonometric functions that results in $\left(\cos\delta_\tx{r},-\sin\delta_\tx{r}\right)$ {(see \cite{TAC2020} which shows how $\delta_\tx{r}$ relates to the power and voltage set-points)}. Finally, by exploiting the angle difference trigonometric identities and the half-angle sine formula \cite[Lemma 1 and Proposition 5]{TAC2020} an approximation of the angle term in \eqref{eq:hac} is obtained
\begin{equation}\label{eq:implementation}
-\sin\left(\dfrac{\delta-\delta_\tx{r}}{2}\right)\approx \dfrac{ \sin\delta_\text{r} \widetilde{\cos\delta}-\widetilde{\sin\delta} \cos\delta_\text{r}}{\sqrt{2\left(1+\widetilde{\cos\delta} \cos\delta_\text{r} + \widetilde{\sin\delta} \sin\delta_\text{r}\right)}}.
\end{equation}  

Respective multiplication of \eqref{eq:implementation} and the dc voltage error with the gains $\gamma_\tx{ac}$ and $\gamma_\tx{dc}$ provides the necessary ingredients for implementing \ac{hac} \eqref{eq:hac}. 
The reader is referred to \cite[Section V]{TAC2020} for further discussions.

\section{Numerical Case Studies}\label{sec:simulation}
\begin{figure}
	\centering
	{\includegraphics[trim=1mm 0mm 0mm 1mm ,clip,width=0.85\columnwidth]{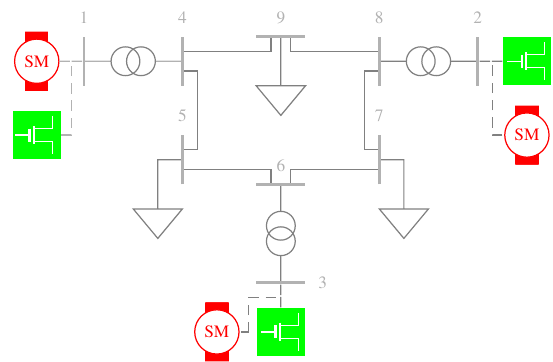}}
	\caption{{IEEE 9-bus test system including SMs and \ac{gfc}s; depending on the choice of the generation technology at nodes 1, 2, and 3 this network model represents: 1) a conventional system (i.e., all-SMs), 2) a SM-dominated system (2SMs-1GFC), 3) a low-inertia system (i.e., 1SM-2GFCs), and 4) a purely converter-based system (i.e., all-GFCs).}}\label{fig:9bus-sys}
\end{figure}

In this section, we employ the standard IEEE 9-bus test system as described in \cite{MOVAH19,TGAKD20,TGA19}; see Figure \ref{fig:9bus-sys} and Table \ref{Table} for the model and control parameters. 
This system represents a transmission grid in which the \ac{sm}s are rated at $100~[\mathrm{MVA}]$. Thus, in order to study the system-level influence of the \ac{gfc}s, it is ideally desired to consider a roughly similar rating for the \ac{gfc}s as that of the \ac{sm}s. 
%
%
{Therefore, we employ the dynamic aggregation technique proposed in \cite{BJJBD20} that provides appropriate scaling laws for deriving the aggregated parameters of a high-power converter system based on the parameters of a smaller module.}%
To this end, we consider a \ac{gfc} module rated at $500~[\mathrm{kVA}]$ that corresponds to a commercially available system; see \cite[Table 1]{TGAKD20} for the parameters. Next, we envision a parallel connection of $200\times500~[\mathrm{kVA}]$ units that is represented with a $100~[\mathrm{MVA}]$ aggregated model with the same structure as in Figure \ref{fig:converter}. The reader is referred to \cite[Figure 2 and Remark 1]{TGAKD20}\cite{CTGAKF19,BJJBD20} for further details and similar model aggregation approaches. In what follows, the system-level performance of the grid-forming hybrid angle control is explored in various case studies. Note that the employed simulation model in MATLAB/Simulink environment is publicly available \cite{TGA19}.
\begin{table}
	\centering
	\caption{Case study model and control parameters \cite{TGA19}.\label{Table}}
	\hspace{-4mm}
	\scalebox{0.85}{
		{\renewcommand{\arraystretch}{1.4}    
			\centering
			\begin{tabular}[]{|c|c||c|c||c|c|}
				\hline\hline
				
				\multicolumn{6}{|c|}{{IEEE 9-bus test system base values}}\\
				\hline
				$S_\text{b}$ & $100$ MVA & $v_\text{b}$ & $230$ kV& $\omega_\text{b}$ & $50$ Hz\\
				\hline
				
				\multicolumn{6}{|c|}{{synchronous machine (SM)}}\\ \hline
				$S_\text{r}$ & $100$ MVA & $v_\text{r}$ & $13.8$ kV& $D$ & {$0$} \\\hline
				$H$ & $3.7$ s & $d_p$ & $1\%$ & $\tau_{\text{g}}$ & $5$ s\\\hline
				
				\multicolumn{6}{|c|}{{single converter module}}\\\hline
				$S_\text{r}$ & $500$ kVA & $G_{{\text{\upshape{dc}}}},C_{{\text{\upshape{dc}}}}$ & $0.83,0.008$ $\Omega^{-1}$,F & $v_{{\text{\upshape{dc,r}}}},v_\text{ll-rms,r}$ & $2.44,1$ kV\\ \hline
				$R$ & $0.001$ $\Omega$ & $L$ & $200$ $\mu$H & $C$ & $300$ $\mu$F\\\hline
				$n$ & $200$ & $\tau_{{\text{\upshape{dc}}}}$ & $50$ ms & $i^{\text{dc}}_{\max}$ & $1.2$ pu\\\hline
				
				\multicolumn{6}{|c|}{{dc and ac voltage controls}}\\\hline
				$k_{{\text{\upshape{dc}}}}$ & $1.6\times10^{3}$&$k_\text{p}$&$0.001$&$k_\text{i}$&$0.5$\\\hline
				
				\multicolumn{6}{|c|}{{droop control, matching control, and HAC}}\\\hline
				$\omega_0$ & $\omega_\tx{b}$& $d_{p-\omega}$ & $1\%$ & $\eta$ & $\omega_0/v_\tx{dc,r}$ \\\hline
				$\gamma_\tx{dc}$ & $0.01\eta$ & $\gamma_\tx{ac}$ & $205$ & $\delta_\tx{r}$ & $0.0238$ \\\hline\hline
				
	\end{tabular}}}
\end{table}

\subsection{Purely Converter-Based System}\label{subsec:all-GFCs}
\begin{figure}
\centering	
{\includegraphics[trim=7.5mm 6mm 14mm 8.5mm ,clip,width=0.88\columnwidth]{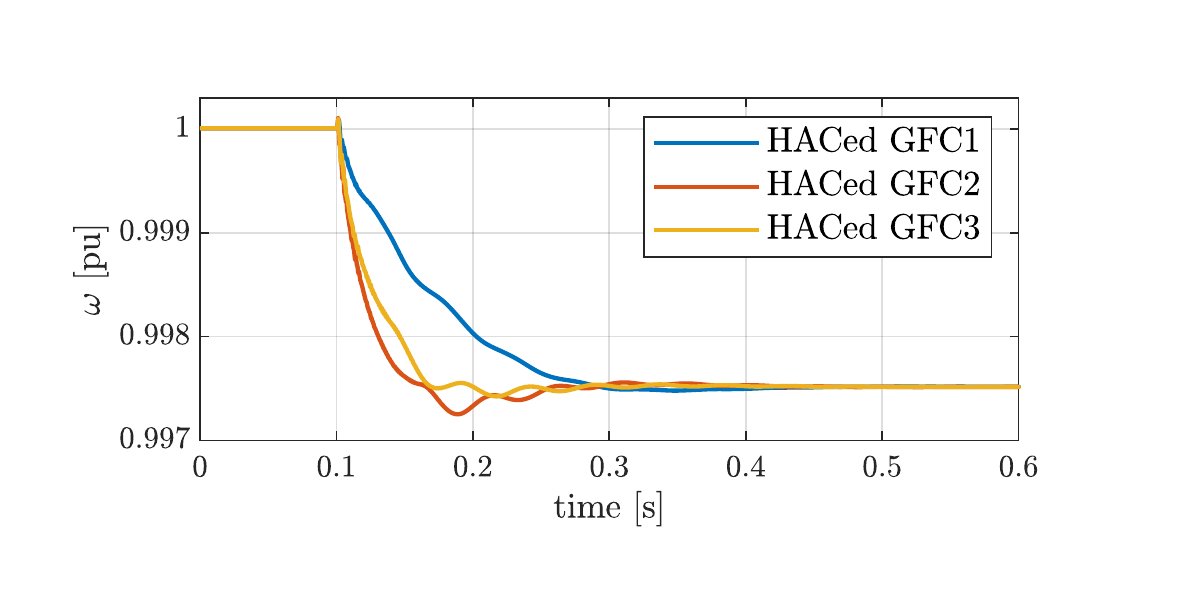}}\vspace{2mm}
{\includegraphics[trim=7.5mm 6mm 14mm 8.5mm ,clip,width=0.88\columnwidth]{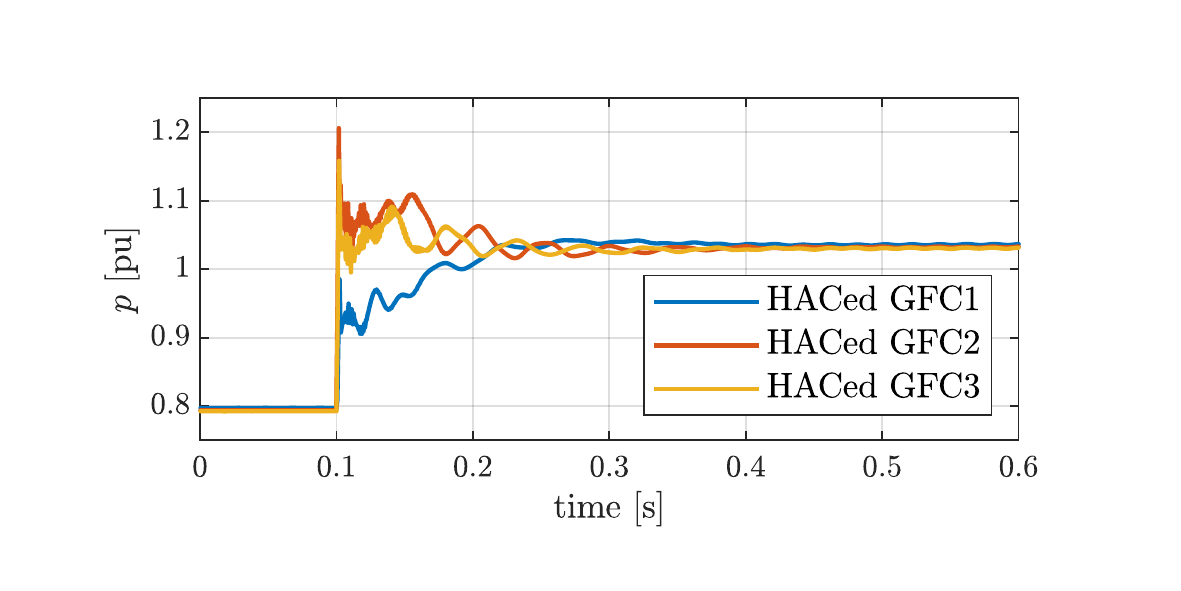}}
\caption{Frequency response of the all-GFCs IEEE 9-bus system configuration under \ac{hac} strategy following a load disturbance (top), the active power time-evolution associated with the GFCs at node 1, 2, and 3 (bottom).}\label{fig:all-GFCs}\vspace{5mm}
{\includegraphics[trim=7.5mm 6mm 14mm 8.5mm ,clip,width=0.88\columnwidth]{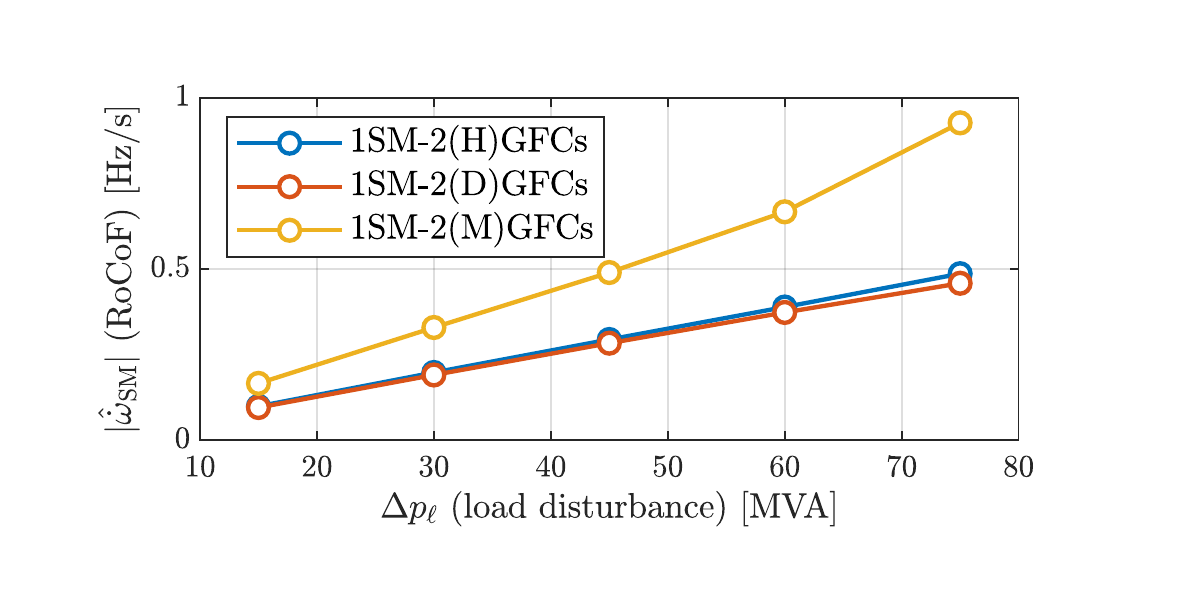}}\vspace{2mm}
{\includegraphics[trim=7.5mm 6mm 14mm 8.5mm ,clip,width=0.88\columnwidth]{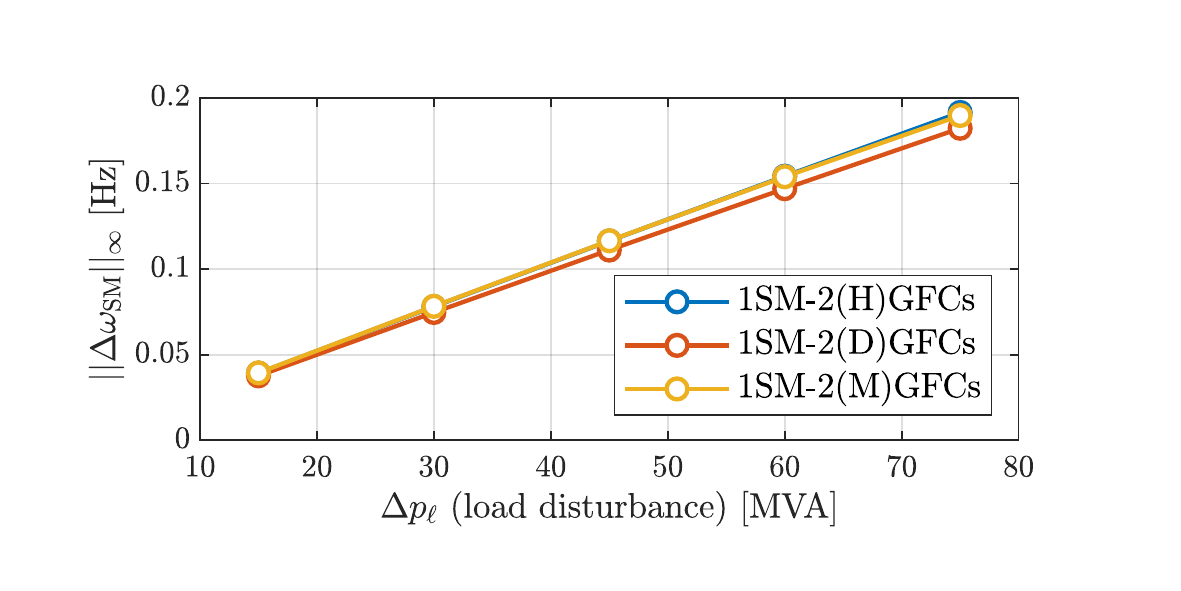}}
\caption{\ac{rocof} evolution of the SM at node 1 in 1SM-2GFCs system configuration under different controls and with respect to variations in network load disturbance (top), maximum frequency deviation (bottom).}\label{fig:metrics}\vspace{5mm}	
{\includegraphics[trim=7.5mm 6mm 14mm 8.5mm ,clip,width=0.88\columnwidth]{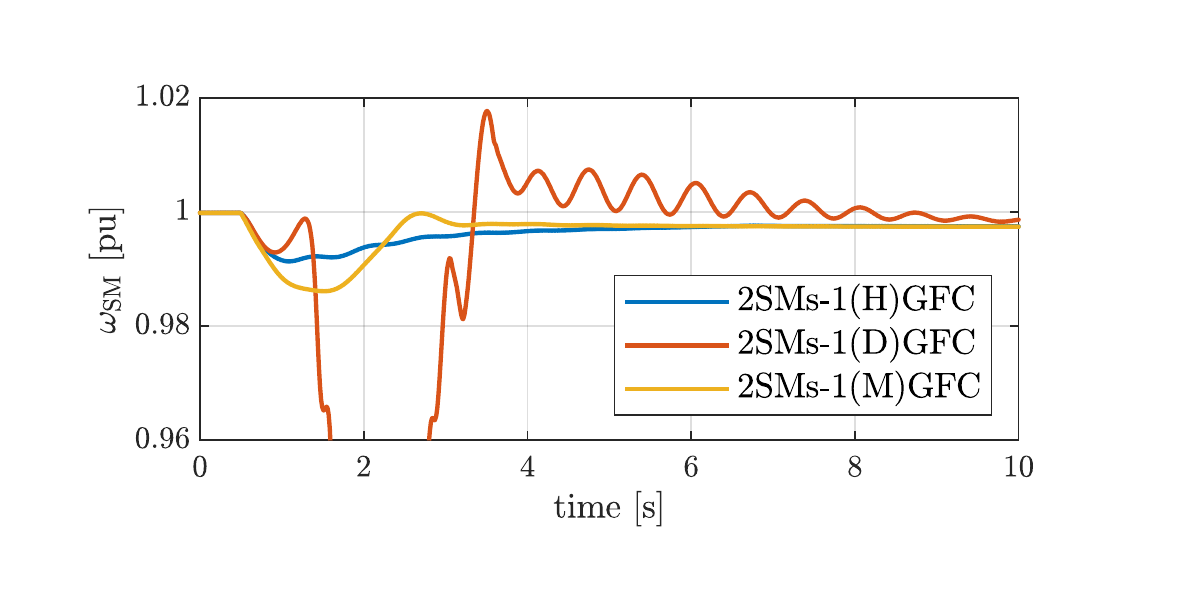}}
\caption{The post-event frequency time-evolution of the SM at node 1 in 2SMs-1GFC system configuration under different controls for the GFC.}\label{fig:instability}
\end{figure}

In this scenario, we remove the \ac{sm}s from the grid model shown in Figure \ref{fig:9bus-sys} and consider three identical \ac{gfc}s at nodes 1-3. We implement the \ac{hac} strategy (as in Figure \ref{fig:hac}) for all the \ac{gfc}s and identically tune the controllers (resulting in equal load-sharing). Figure \ref{fig:all-GFCs} illustrates the frequency and active power evolution at the generation nodes following a {$0.75~[\mathrm{pu}]$} load-disturbance at node 7. Note that we do not implement an explicit frequency measurement mechanism and rather observe the internal frequency of the GFCs i.e., \eqref{eq:hac}. 

Figure \ref{fig:all-GFCs} firstly verifies the grid-forming nature of \ac{hac} technique i.e., the autonomous operation of converters {without relying on a reference frequency provided by an external source (as in the case of grid-following devices)}. Second, note that all GFCs reach the post-contingency equilibrium in approximately $200~[\mathrm{ms}]$. This fast response timescale is crucial for a satisfactory grid-forming performance \cite{MDHHV18}. Observe that the units that are {electrically} closer to the disturbance location i.e., GFCs at node 2 and 3 react on a slightly faster timescale. Last, Figure \ref{fig:all-GFCs} confirms the equal load-sharing of the \ac{gfc}s. See \cite{TGAKD20,CTGAKF19} for similar observations concerning the behavior of other grid-forming schemes in an all-GFCs network. 

\subsection{{The Influence on the Frequency Performance Metrics}}\label{subsec:metrics}
In what follows, we explore the influence of grid-forming architectures (presented in Section \ref{sec:control}) on the frequency stability of low-inertia configuration associated with the grid model shown in Figure \ref{fig:9bus-sys} i.e., with \ac{sm} at node 1 and \ac{gfc}s at nodes 2 and 3. More precisely, we 1) implement identical droop controllers for both GFCs, 2) implement matching control for the GFCs, and finally 3) consider \ac{hac} for the converters. Subsequently, for all the 1SM-2GFCs pairs under different controls, we apply five load-disturbances at note 7 i.e., $\Delta p_{\ell,k}=15k~[\mathrm{MVA}]$ where $k=1,\ldots,5$. Next, we observe the frequency dynamics via the angular frequency of the \ac{sm} at node 1; see \cite{TGAKD20,CTGAKF19} for similar approaches. We evaluate the standard frequency performance metrics, namely: 1) maximum frequency deviation, and 2) \ac{rocof} that are (with slight abuse of notation) defined by \cite{poolla_placement_2018,TGAKD20}
\begin{subequations}\label{eqs:metrics}
\begin{align}
\n{\Delta\omega_\tx{SM}}_\infty&:=\max\left(\omega_0-\omega_\tx{SM}(t)\right),\\
\big|\hat{{\dot{\omega}}}_\tx{SM}\big|&\coloneqq\bigg|\frac{\omega_\tx{SM}(t_0+\Delta t)-\omega_\tx{SM}(t_0)}{\Delta t}\bigg|,
\end{align}	
\end{subequations}     
where $\omega_\tx{SM}$ and $\Delta t$ respectively denote the \ac{sm} frequency and the \ac{rocof} approximation time horizon. We consider $\Delta t=150~[\mathrm{ms}]$ that allows to observe the influence of fast GFCs dynamics on the system frequency; the reader is referred to \cite{CTGAKF19} for a detailed discussion on the choice of RoCoF window.

Figure \ref{fig:metrics} highlights the evolution of metrics \eqref{eqs:metrics} in 1SM-2GFCs system, under different grid-forming controls, and with respect to the disturbance variation. Concerning the \ac{rocof}, the techniques incorporating an ac feedback i.e., droop control \eqref{eq:droop} and HAC \eqref{eq:hac} result in better performance compared to the purely dc-based matching control. This is underpinned by their fast disturbance sensing and frequency modification feature. In contrast, matching control \eqref{eq:matching} only reacts when the network disturbance is propagated to the converter dc dynamics. 

Figure \ref{fig:metrics} confirms that \ac{hac} inherits the advantage of the other ac-based controls e.g., droop control; see Remark \ref{rem:ac} and \cite{TGAKD20} for further details. On the other hand, maximum frequency deviation performance is almost identical when employing different controls. The reader is referred to \cite{TGAKD20} for a detailed discussion on this aspect. Although Figure \ref{fig:metrics} uncovers a performance variation that depends on the grid-forming controls in 1SM-2GFCs setup, all controls outperform the all-SMs configuration by a clear margin; e.g., see \cite[Figures 11 and 12]{TGAKD20}. The positive influence of the grid-forming controls on frequency stability is also reported in \cite{CTGAKF19,MOVAH19}.  

\subsection{Synchronous Machine Dominated System}
In this test case, we consider a SM-dominated generation profile by including identical SMs at nodes 1 and 3, and a GFC at node 2 in Figure \ref{fig:9bus-sys} leading to a 2SMs-1GFC system. We consider the same load-disturbance as in Subsection \ref{subsec:all-GFCs} and enforce equal load-sharing for all units. Figure \ref{fig:instability} shows the frequency of the SM at node 1 when the GFC is under different controls. The network base load and disturbance are chosen such that the desired post-disturbance GFC's active power does not result in exceeding the converter dc source current limit; see Figure \ref{fig:converter}. However, due to the presence of slowly reacting SMs (because of large turbine time constants), the GFC dominantly supports the load during the transient that results in violating the dc limit. As it is also reported in \cite{SC21,GRL20,TGAKD20}, under active dc current constraint, droop control exhibits instability by depleting the dc-link energy. This aggressive behavior is because droop control \eqref{eq:droop} is agnostic to the dc dynamics. Subsequently, the converter instability propagates through the network an destabilizes the SMs. This adverse interaction can be counteracted by augmenting droop control with an ac current limiting mechanism e.g., threshold virtual impedance \cite{GD19}. However, in our experience the limiting performance is fragile and depends on various factors e.g., the disturbance size and location. 

In contrast, as it can be observed in Figure \ref{fig:instability}, the controls that exploit dc feedback in their angle dynamics e.g., matching control and \ac{hac} exhibit robustness with respect to the dc current constraint. This behavior is numerically and formally explored in \cite{TGAKD20,SC21,GRL20}. In a nutshell, it is underpinned by an inherent mode-switching behavior that transforms the grid-forming control to a following one and allows injecting constant current while dc source is saturated. The reader is referred to \cite{TGAKD20} for a detailed discussion. Interestingly, \ac{hac} results in better frequency deviation performance compared to matching control. This is because \ac{hac} provides global stability certificates \cite{TAC2020,TAD20} and includes a hybrid ac/dc structure that strikes a balance between robustness and performance.           
\begin{remark}\textup{(Adverse timescales interaction)}\\
In case studies presented in \cite{TGAKD20,SC21,GRL20}, the instability of the purely ac-based grid-forming techniques strongly depends on the disturbance magnitude. However, Figure \ref{fig:instability} provides an alternative insight. More precisely, the base load and disturbance scenario specification does not derive the converter beyond its limit in all-GFCs or 1SM-2GFCs systems. However, in SM-dominated system the adverse interplay of the fast GFC and slow SMs results in excessive transient current injection by the converter. Hence, aside the disturbance characteristics, the penetration level of the converter-based generation is an influencing factor for this instability mechanism.	
\end{remark}
\subsection{{The Influence of GFC parameters on Frequency Behavior}}
Several works e.g., \cite{RPG18,d2013virtual,JRSD17} highlight the necessity and benefits of low-pass filtering the ac measurements that are exploited in synthesizing converter control. In the following, we show how an appropriate tuning of the converters' low-pass filters improves the frequency stability and can be perceived as a decentralized frequency shaping strategy.

\begin{figure}
	\centering
	{\includegraphics[trim=7.5mm 6mm 14mm 8.5mm ,clip,width=0.88\columnwidth]{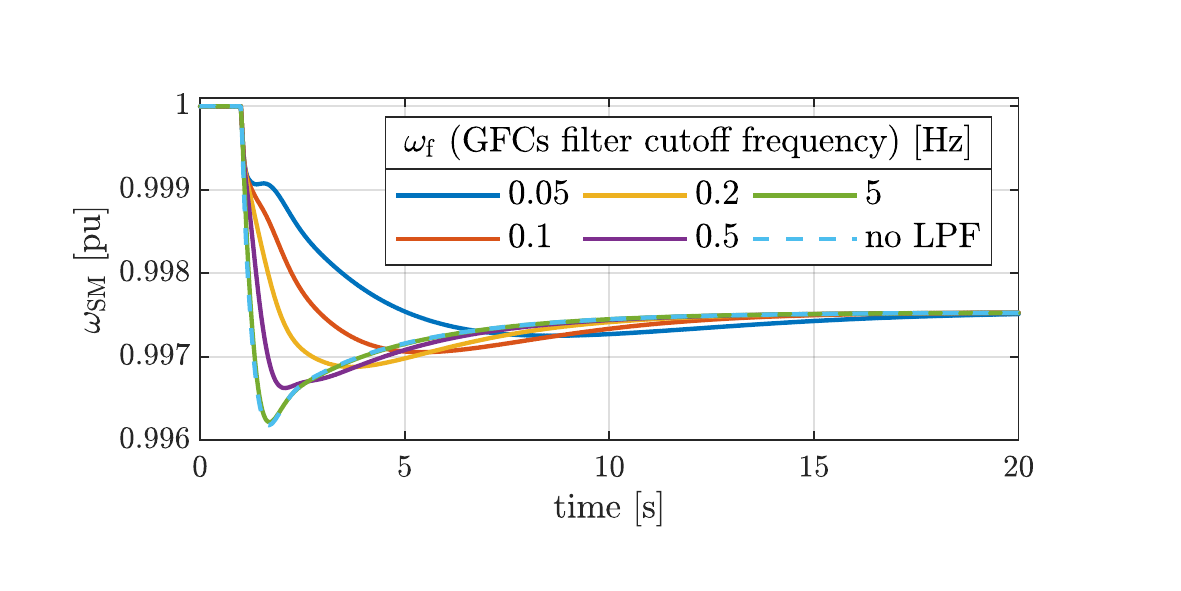}}
	\caption{Frequency evolution of SM at node 1 when the GFCs are controlled by HAC strategy with 5 different LPF cutoff frequency and without a LPF.  }\label{fig:LPF}\vspace{5mm}
	{\includegraphics[trim=7.5mm 6mm 14mm 8.5mm ,clip,width=0.88\columnwidth]{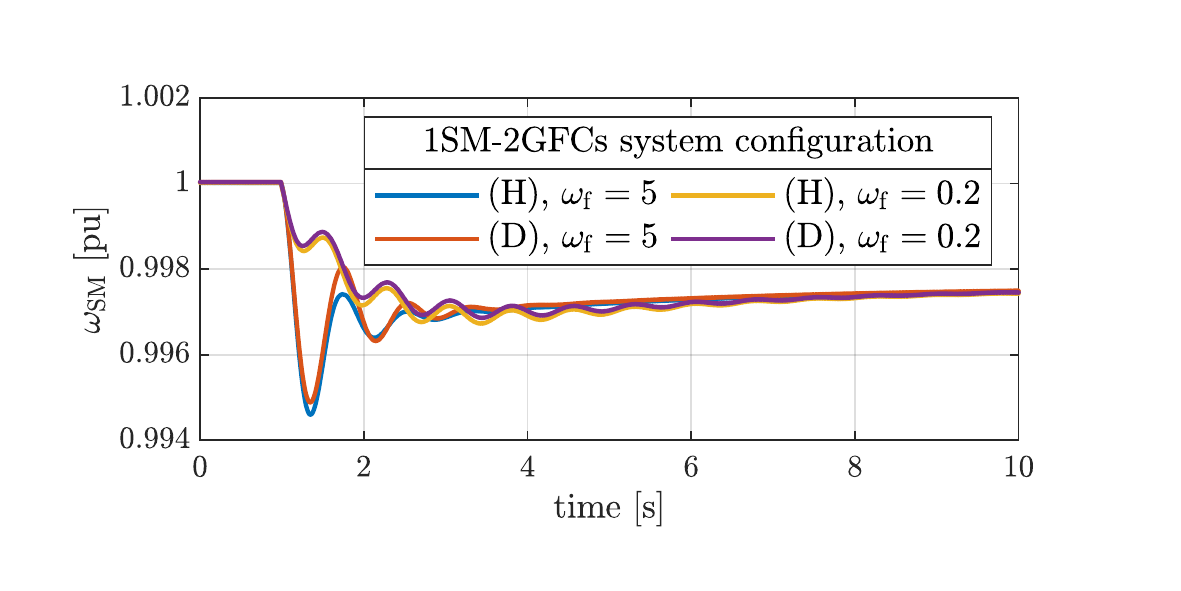}}
	\caption{The PSS-free response of the SM at node 1 when the GFCs are controlled by either HAC or droop control while considering two different LPF cutoff frequencies.}\label{fig:PSS-free}
\end{figure}

To begin with, we consider the 1SM-2GFCs configuration as in Subsection \ref{subsec:metrics}. Next, we equip the GFCs at nodes 2 and 3 with HAC and consider different cutoff frequencies for the underlying low-pass filters; see Figure \ref{fig:hac}. Figure \ref{fig:LPF} illustrates the frequency response of the SM with respect to the converters' cutoff frequency variation. Note that this observation suggests that an appropriate tuning of the LPFs can globally reshape the frequency response across the system. This is due to the fact that the LPFs with an appropriate tuning provide a certain amount of virtual inertia and thus enhance the SM frequency response. The reader is referred to \cite{d2013virtual} for detailed discussion. We remark that the positive influence of $\omega_\tx{f}$-variation on frequency stability must be compared against the enhancement that is achieved via changing the converters' controls. In particular, the maximum frequency deviation metric is significantly reduced in Figure \ref{fig:LPF} as $\omega_\tx{f}$ is reduced, in contrast, for a fixed $\omega_\tx{f}$ the improvement due to changing grid-forming strategies is minimal; see Figure \ref{fig:metrics}. Last, our numerical investigations confirms an almost identical behavior as when the GFCs are controlled by droop control.

\begin{figure*}
	\centering
	\includegraphics[width=\linewidth]{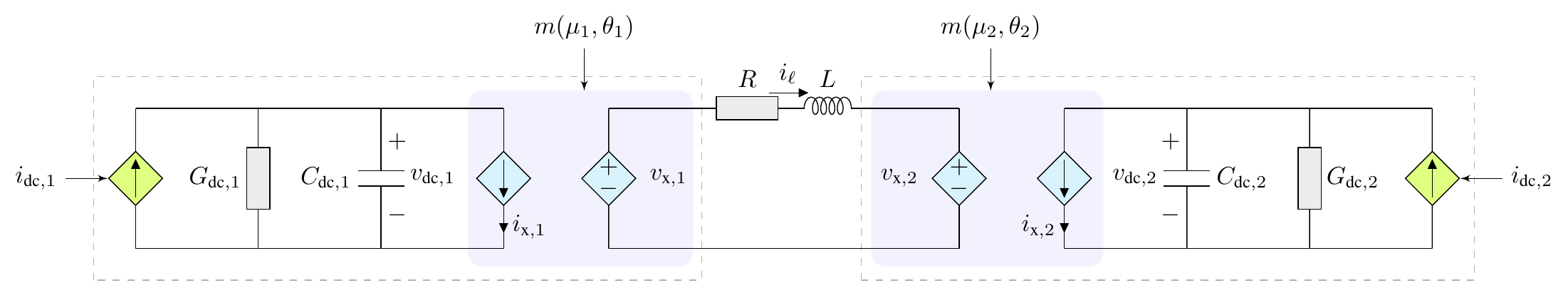}
	\caption{The simplified two-converter system configuration where $i_{\tx{x},j}:=m(\mu_j,\theta_j)^\top i_\ell$ and $v_{\tx{x},j}:=m(\mu_j,\theta_j) v_{\tx{dc},j}$ for $j=1$ and $2$.}\label{fig:two-converter}
\end{figure*}
\subsection{The PSS-Free Behavior of the Low-Inertia System}
It has been recently uncovered that the \ac{pss}s of the SMs in a low-inertia system might adversely interact with certain timescales of the converter systems \cite{MOVAH19}. Moreover, it is highlighted that an adaptive PSS might be a crucial element in ensuring system stability in transitioning to 100$\%$ converter-based generation \cite{MMF20,CTGAKF19,RK20}. However, the online/offline modification of the well-established PSS control architectures might be an expensive, challenging, and practically infeasible task. {In this subsection, we investigate if the LPF retuning strategy can possibly allow for a \ac{pss}-free operation of the SMs in a low-inertia network configuration.}  
%

To this end, we consider the 1SM-2GFCs configuration and remove the PSS form the SM model at node 1. Similar to previous cases, we separately consider droop control and HAC for the GFCs. Next, we consider two different cutoff frequencies for the underlying LPFs. Figure \ref{fig:PSS-free} depicts the frequency response of the SM. It can be seen that under $\omega_\tx{f}=5~[\mathrm{Hz}]$, HAC results in slightly less oscillations compared to droop control. However, when selecting $\omega_\tx{f}=0.2~[\mathrm{Hz}]$, the SM frequency nadir (that might potentially trip the low-frequency protection mechanisms) is almost removed. Although, under $\omega_\tx{f}=0.2~[\mathrm{Hz}]$ the low-frequency oscillation are still present but the appropriate LPF tuning reshapes the envelope on the frequency response and constrains the frequency oscillations within the stability margin. Last, the enhanced response is achieved regardless of the control strategies of the GFCs. {We remark that this insight is a result of a preliminary exploration and requires an in-depth analysis.}      
\section{Conclusions and Outlook}
In this paper, the system-level performance and robustness of grid-forming \ac{hac} is explored. We verified that the multi-variable dc-ac \ac{hac} inherits the enhanced performance and robustness of ac-based and dc-based controls, respectively. Furthermore, we highlighted that how retuning the low-pass filters for grid-forming controls can significantly enhance the frequency performance across the system. Last, an alternative augmented controller is introduced and a simplified stability analysis for a two-converter system is presented. Our future work includes: 1) implementation of \ac{hac} for power converters in more complex power grid models, 2) power hardware validation of \ac{hac} concept, and 3) detailed stability analysis multi-converter system under \ac{hac}.
\appendix
\subsection{Dynamic Inverse Droop Control Augmentation}
The \ac{hac} implementation described in the Section \ref{sec:control}, explicitly relies on the foreknowledge of the relative angle reference $\delta_\tx{r}$. The reader is referred to \cite[Proposition 6]{TAC2020} on deriving $\delta_\tx{r}$ based on the power and voltage set-points. However, this techniques relies on the system parameters and might exhibit robustness issues. In what follows, we introduce a complementary high-level feedback control mechanism that achieves the same objective with enhanced robustness.

To begin with we introduce an integral control that defines $\delta_\tx{r}$ in relation to the active power mismatch i.e., 
\begin{equation}\label{eq:p-delta}
\delta_\tx{r}=\int_{0}^{t}\kappa_{p\delta}(p_\tx{r}-p),
\end{equation}
where $\kappa_{p\delta}$ denotes the integrator gain. It is worth mentioning that augmenting the \ac{hac} with \eqref{eq:p-delta} transforms the grid-forming \ac{hac} to a phase locked loop (PLL)-free grid-following controller that achieves robust active power reference tracking. This is due to the fact that the integrator \eqref{eq:p-delta} disables the natural $p-\omega$ droop behavior of the \ac{hac} \cite[Proposition 7]{TAC2020}. 

Although robust reference tracking and disturbance rejection might be desirable for certain applications, drooping behavior that enables load-sharing between generation units in transmission grids is vital. Thus, we combine \eqref{eq:p-delta} with an inverse $\omega-p$ droop control (reminiscent of the speed droop control of the \ac{sm} e.g., see \cite[Equation 5]{TGAKD20}) that is  
\begin{equation}\label{eq:inverse droop}
p_\tx{r}:=p^\star+d_{\omega-p}(\omega_\tx{c}-\omega_0),
\end{equation}
where $p^\star$ denotes the power reference at nominal frequency, $d_{\omega-p}$ is the inverse droop gain, and $\omega_\tx{c}$ is the internal feedback given by \ac{hac} \eqref{eq:hac}. We emphasize that the combination of \eqref{eq:hac} with the dynamic inverse droop control \eqref{eq:p-delta} and \eqref{eq:inverse droop} is a heuristic approach, and requires a separate detailed stability analysis as in \cite{TAC2020}. However, our numerical investigations confirms that the cascade control structure achieves the aforementioned control specification. Furthermore, as in standard cascaded control systems (e.g., \cite{SGCD19}) the controls \eqref{eq:hac}, \eqref{eq:p-delta}, and \eqref{eq:inverse droop} must be tuned in harmony while respecting the required timescales separation between the nested loops. Last, we remark that since the main focus of this work is to explore the behavior of the standard \ac{hac} \eqref{eq:hac}, the presented case studies in the Section \ref{sec:simulation} only incorporates \eqref{eq:hac}.  

\subsection{Stability Analysis of the Interconnected Converters}
The stability analysis of the \ac{hac} in \cite{TAC2020,TAD20} is centered around a model configuration that includes a single converter connected to an infinite bus or a dynamic center of inertia (CoI) grid model. The converter-CoI model abstractly represents the mixed converter-SM configurations of the IEEE 9-bus system shown in Figure \ref{fig:9bus-sys}. However, the stability and synchronization analysis of a multi-converter setup requires a separate investigation and is presented as it follows.

To begin with, we consider a simplified two-converter model as in Figure \ref{fig:two-converter}; see \cite{JD20} for a similar model configuration. Without loss of generality we assume the line current is flowing out of the converter 1 and into the converter 2. It is noteworthy that the configuration in Figure \ref{fig:two-converter} represents two converters with RL output filters that are connected through a RL line section. Thus, the filter elements and the line section are merged into a single RL elements.

We equip both converters with \ac{hac} \eqref{eq:hac} and thus the overall dynamical model is given by
\begin{subequations}\label{eqs:sys}
	\begin{align}
	&\dot{\theta}_1=\omega_0+\gamma_{\tx{dc},1}(v_{\tx{dc},1}-v_{\tx{dc,r}1})-\gamma_{\tx{ac},1}\sin\left(\dfrac{\delta-\delta_r}{2}\right),\label{eq:sys1}
	\\
	&\dot{v}_{\tx{dc},1}=C_{\tx{dc},1}^{-1}\left(i_{\tx{dc},1}-G_{\tx{dc},1}v_{\tx{dc},1}-m(\mu_1,\theta_1)^\top i_\ell\right),\label{eq:sys2}
	\\
	&\dot{\theta}_2=\omega_0+\gamma_{\tx{dc},2}(v_{\tx{dc},2}-v_{\tx{dc,r}2})+\gamma_{\tx{ac},2}\sin\left(\dfrac{\delta-\delta_r}{2}\right),\label{eq:sys3}
	\\
	&\dot{v}_{\tx{dc},2}=C_{\tx{dc},2}^{-1}\left(i_{\tx{dc},2}-G_{\tx{dc},2}v_{\tx{dc},2}+m(\mu_2,\theta_2)^\top i_\ell\right),\label{eq:sys4}
	\end{align}
\end{subequations}
where $\delta:=\theta_1-\theta_2$. In addition, the dc source current and the modulation signal control input pairs are defined as
\begin{subequations}\label{eqs:controls}
	\begin{align}
	i_{\tx{dc},1}&:=-\kappa_{\tx{dc},1}\left(v_{\tx{dc},1}-v_{\tx{dc,r}1}\right),
	\\
	m(\mu_1,\theta_1)&:=\mu_1 \left(\cos\theta_1,\sin\theta_1\right)^\top,
	\\
	i_{\tx{dc},2}&:=-\kappa_{\tx{dc},2}\left(v_{\tx{dc},2}-v_{\tx{dc,r}2}\right),
	\\
	m(\mu_2,\theta_2)&:=\mu_2 \left(\cos\theta_2,\sin\theta_2\right)^\top.
	\end{align}
\end{subequations}
Furthermore, for simplicity of exposition we assume quasi-steady-state line dynamics, thus the current in stationary $\alpha\beta$-coordinates is given by
\begin{equation}\label{eq:line current}
i_\ell:=R^{-1}\left(v_{\tx{dc},1}m(\mu_1,\theta_1) - v_{\tx{dc},2}m(\mu_2,\theta_2)\right).
\end{equation}
We subsequently 1) combine the absolute angle dynamics \eqref{eq:sys1} and \eqref{eq:sys3} into relative angle dynamics, 2) replace $i_\ell$ in \eqref{eq:sys2} and \eqref{eq:sys4} with the expression in \eqref{eq:line current}, and 3)  incorporate the controls \eqref{eqs:controls} that results in the overall closed-loop dynamics:
\begin{subequations}\label{eqs:closed-loop}
	\begin{align}
	\dot{v}_{\tx{dc},1}=&- C_{\tx{dc},1}^{-1} \Big(\kappa_{\tx{dc},1} \left( v_{\tx{dc},1} - v_{\tx{dc,r}1} \right) + \left(G_{\tx{dc,1}}+R^{-1}\mu_1^2\right)v_{\tx{dc},1}\nonumber
	\\
	& -  R^{-1}\mu_1\mu_2\cos(\delta)v_{\tx{dc},2}\Big),\label{eq:sys2modified}
	\\
	\dot{v}_{\tx{dc},2}=&- C_{\tx{dc},2}^{-1} \Big(\kappa_{\tx{dc},2} \left( v_{\tx{dc},2} - v_{\tx{dc,r}2} \right) + \left(G_{\tx{dc},2}+R^{-1}\mu_2^2\right)v_{\tx{dc},2}\nonumber
	\\
	& -  R^{-1}\mu_1\mu_2\cos(\delta)v_{\tx{dc},1}\Big),\label{eq:sys4modified}
	\\
	\dot{\delta}=&+\gamma_{\tx{dc},1}(v_{\tx{dc},1}-v_{\tx{dc,r}1})-\gamma_{\tx{dc},2}(v_{\tx{dc},2}-v_{\tx{dc,r}2})
	\nonumber
	\\
	&-(\gamma_{\tx{ac},1}+\gamma_{\tx{ac},2})\sin\left(\dfrac{\delta-\delta_r}{2}\right).
	\end{align}
\end{subequations}
Assume that there exists an equilibrium and it coincides with the desired control references (similar to \cite[Theorem 1]{TAC2020}) i.e., 
\begin{equation}
\left(\delta^\star,v_{\tx{dc},1}^\star,v_{\tx{dc},2}^\star\right):=\left(\delta_\tx{r},v_{\tx{dc,r}1},v_{\tx{dc,r}1}\right).
\end{equation}
Next, we associate an energy function with the closed-loop dynamics \eqref{eqs:closed-loop}, that is given by
\begin{align}
\mathcal{V}:=&\frac{1}{2}\left(C_{\tx{dc},1}\left(v_{\tx{dc},1}-v_{\tx{dc,r}1}\right)^2 + C_{\tx{dc},2}\left(v_{\tx{dc},2}-v_{\tx{dc,r}2}\right)^2 \right)
\nonumber
\\
& + 2\left(1-\cos\left(\dfrac{{\delta-\delta_\tx{r}}}{2}\right)\right).
\end{align}
Following the analysis recipe as in \cite[Theorem 2]{TAC2020}, it is possible to show that there exists positive lower bounds $\gamma_{\tx{ac},\min}$ and $\kappa_{\tx{dc},\min}$  such that if $\gamma_{\tx{ac},1}+\gamma_{\tx{ac},2}>\gamma_{\tx{ac},\min}$ and $\kappa_{\tx{dc},1},\kappa_{\tx{dc},2}>\kappa_{\tx{dc},\min}$ then the energy function $\mathcal{V}\to 0$ for (almost) all the solutions of \eqref{eqs:closed-loop}. Note that when $\mathcal{V}\to 0$ it yields that $\left(\delta,v_{\tx{dc},1},v_{\tx{dc},2}\right)\to\left(\delta_\tx{r},v_{\tx{dc,r}1},v_{\tx{dc,r}1}\right)$. This result verifies the asymptotic convergence of converter dynamics to the desired dc voltage and relative angle references. 
\bibliographystyle{IEEEtran}
\bibliography{IEEEabrv,Ref_PSCC}
\end{document}